\documentclass[12pt]{article}

\usepackage[english]{babel}
\usepackage[T1]{fontenc}
\usepackage{amsmath}
\usepackage{amsfonts}
\usepackage{amssymb}
\usepackage[utf8]{inputenc}
\usepackage{graphicx,color,url}
\usepackage[margin=2.5cm]{geometry}

\def\no{\noindent}
\newtheorem{lem}{Lemma}

\newtheorem{theo}{Theorem}
\newtheorem{rem}{Remark}
\def\RR{\mathbb{R}}
\def\ZZ{\mathbb{Z}}
\def\bz{\mathbf{z}}

\def\bq{\mathbf{q}}
\def\bp{\mathbf{p}}

\def\bu{\mathbf{u}}
\def\bb{\mathbf{b}}
\def\bc{\mathbf{c}}

\def\bpsi{{\boldsymbol{\psi}}}

\def\bom{{\boldsymbol{\omega}}}

\def\calH{\mathcal{H}}

\def\pmatrix{\left(\begin{array}}
\def\endpmatrix{\end{array}\right)}
\def\dd{\mathrm{d}}
\def\sech{\mathrm{\,sech}}
\def\atan{\mathrm{\,atan}}

\def\dbq{\dot{\bq}}
\def\dbp{\dot{\bp}}

\def\ddgam{\ddot{\gamma}}
\def\ddeta{\ddot{\eta}}

\def\dgam{\dot{\gamma}}
\def\deta{\dot{\eta}}

\def\proof{\underline{Proof}\quad}
\def\QED{\mbox{\,$\Box{~}$}}
\def\P{{\cal P}}
\def\hX{\hat{X}}
\def\hu{{\hat{u}}}
\def\bu{{\bar{u}}}
\def\bau{{\bar{u}^{(1)}}}
\def\bbu{{\bar{u}^{(2)}}}
\def\baf{{\bar{f}}}

\def\bpsi{{\bar{\psi}}}

\def\aa{\alpha}
\def\hq{\hat{q}}
\def\hp{\hat{p}}
\def\hs{\hat{s}}
\def\hc{\hat{c}}

\title{Energy conserving methods for Hamiltonian PDEs based on spectral space decomposition}

\author{L.\,Brugnano$^a$\quad G.\,Frasca Caccia$^a$\quad F.\,Iavernaro$^b$\\[.5cm]
\small
$^a$\,Dipartimento di Matematica e Informatica ``U.\,Dini'', Universit\`a di Firenze, Italy\\
\small
$^b$\,Dipartimento di Matematica, Universit\`a di Bari, Italy}

\begin{document}

\maketitle

\abstract{In this paper we discuss energy conservation issues related to the  numerical solution of the nonlinear wave equation, when a Fourier expansion is considered for the space discretization. The obtained semi-discrete problem is then solved in time by means of energy-conserving Runge-Kutta methods in the HBVMs class.

\medskip
\no{\bf Keywords:} nonlinear wave equation; Hamiltonian PDEs; Fourier expansion; energy-conserving methods; Hamiltonian Boundary Value Methods; HBVMs.

\medskip
\no{\bf AMS:} 65P10, 65L05, 65M20.}

\section{Introduction}\label{intro}

In this paper we discuss energy-conservation issues concerning the nonlinear wave equation, even though the arguments can be extended to different types of Hamiltonian PDEs.  For sake of simplicity,  we shall consider the 1D case,
\begin{eqnarray}\nonumber 
u_{tt}(x,t)&=&\aa^2u_{xx}(x,t)-f'(u(x,t)), \qquad (x,t)\in (0,1)\times(0,\infty),\\ \label{wave}
u(x,0)&=&\psi_0(x),\\ \nonumber
u_t(x,0)&=&\psi_1(x), \qquad x\in(0,1),
\end{eqnarray}
coupled with suitable boundary conditions. As usual, subscripts denote  partial derivatives.
In (\ref{wave}), the functions $f$, $\psi_0$ and $\psi_1$ are supposed to be suitably regular, so they define a regular solution $u(x,t)$ ($f'$ denotes the derivative of $f$).   The problem is completed by assigning suitable boundary conditions which we shall, at first, assume to be periodic,
\begin{equation}\label{perbc}
u(0,t)=u(1,t), \qquad t>0.
\end{equation}
Later on, we shall also consider the case of Dirichlet boundary conditions
\begin{equation}\label{diribc}
u(0,t) = \phi_0(t), \qquad u(1,t) = \phi_1(t), \qquad t>0,
\end{equation}
and Neumann boundary conditions
\begin{equation}\label{neubc}
u_x(0,t) = \phi_0(t), \qquad u_x(1,t) = \phi_1(t), \qquad t>0,
\end{equation}
with $\phi_0(t)$ and $\phi_1(t)$ suitably regular.  In all cases, all the functions are assumed to satisfy suitable compatibility conditions, depending on the considered set of boundary conditions. 
\begin{rem}
It is worth mentioning that a problem defined on a generic interval $[a,b]$, could be always transformed to the form (\ref{wave}), by means of a linear transformation of the $x$ variable. In such a case, the leading coefficient $\aa$ in (\ref{wave}) changes accordingly (i.e., it becomes $(b-a)^{-1}\aa$).
\end{rem}

By setting 
\begin{equation}\label{v}
v = u_t,
\end{equation}
and defining the functional
\begin{equation}\label{H1}
\calH[u,v](t)=\int_0^1 \left[\frac{1}{2}v^2(x,t)+\frac{1}{2}\aa^2u_x^2(x,t)+f(u(x,t))\right]\dd x \equiv \int_0^1 E(x,t)\,\dd x,
\end{equation}
we can rewrite (\ref{wave}) as the infinite-dimensional Hamiltonian system (for sake of brevity, we neglect the arguments of the functions $u$ and $v$)
\begin{equation}\label{Hwave}
\bz_t=J\frac{\delta \calH}{\delta \bz}, 
\end{equation}
where
\begin{equation}\label{Hwave1}
J= \pmatrix  {cc}
0 & 1 \\
-1 & 0
\endpmatrix,
\qquad
\bz=\pmatrix{c} u\\v\endpmatrix,
\end{equation}
and 
\begin{equation}\label{Hwave2}
\frac{\delta \calH}{\delta \bz}=\left(\frac{\delta \calH}{\delta u},\frac{\delta \calH}{\delta v}\right)^\top 
\end{equation}
is the functional derivative of $\calH$ \cite{BFI14}. Indeed, one proves that  (\ref{Hwave})--(\ref{Hwave2}) are equivalent to (\ref{wave}): 
$$\bz_t=\pmatrix{c}u_t\\ v_t\endpmatrix=J\frac{\delta \calH}{\delta \bz}=\pmatrix{c}\frac{\delta \calH}{\delta v}
\\[3mm]-\frac{\delta \calH}{\delta u}\endpmatrix=\pmatrix{c}v\\ \aa^2 u_{xx}- f'(u)\endpmatrix,$$
or
\begin{eqnarray*}\nonumber
u_t (x,t)&=& v(x,t),   \qquad (x,t)\in (0,1)\times(0,\infty),\\ \label{wave1} 
\\[-2mm] \nonumber
v_t(x,t)&=&\aa^2 u_{xx}(x,t)-f'(u(x,t)),
\end{eqnarray*}
that is, the first-order formulation of the first equation in (\ref{wave}). 

In the last decades there has been a growing interest in the numerical treatment of Hamiltonian PDEs arising in many application fields, such as meteorology and weather prediction, quantum mechanics and nonlinear optics \cite{BrRe01}. For this purpose, different approaches have been developed such as multisymplectic methods \cite{BrRe01, FMR06, IsSc04}, splitting methods \cite{F12}, and semi-discretizations by means of the method of lines (MOL).

When the MOL approach is used, the spatial derivatives can be approximated by finite differences (see for example \cite{BFI14}), but a different technique is that of solving the boundary value problem in space by means of spectral methods \cite{CHQZ88, B01, FW78, WMGSS91, LXW11, Shen94, Shen95}. In both cases one can integrate the resulting system in time through suitable standard integrators, though the use of symplectic and/or symmetric methods 
is preferable (see, e.g., \cite{CHL08}).

In particular, in this paper we use Fourier-Galerkin spectral methods for the space semi-discretization of (\ref{wave}). Galerkin methods require to expand the solution of the problem along a basis in which every component satisfies the associated boundary conditions.

When the problem at hand is coupled with the periodic boundary conditions (\ref{perbc}), a trigonometric basis is usually preferred (see for example \cite{CHQZ88, FW78, WMGSS91}). 
If one has to deal with homogeneous boundary conditions, a basis composed by an appropriate combination of Jacobi polynomials can be also considered, as done for example in \cite{LXW11, Shen94}, where a combination of Legendre polynomials is used, or in \cite{Shen95}, where Chebyshev polynomials are employed. 

In the case of general inhomogeneous boundary conditions, a Galerkin method can still be used by considering a suitable  {\em boundary adapted} basis  \cite{CHQZ88}. Alternatively (see, e.g.,  \cite{B01}),  one may transform the problem at hand into an equivalent one having homogeneous boundary conditions. The solution of this equivalent problem can then be expanded along a suitable trigonometric basis. This latter approach is considered in the sequel. 

An important feature that one could be interested to numerically reproduce, is that the variation of the energy density, integrated over an interval, depends only on the net flux through its endpoints. In particular, if there is no net flux (as in the case, for example,  of periodic boundary condition), then the integrated energy density is exactly conserved, meaning that it remains constant over time.

In this paper we show that the use of energy-conserving methods in time, assures a precise reproduction of the above mentioned conservation law of the semi-discrete model obtained by means of a Fourier-Galerkin method in space.
In particular, we shall here consider methods in the class of {\em Hamiltonian Boundary Value Methods (HBVMs)}, recently introduced for the numerical solution of Hamiltonian problems \cite{BIT09,BIT09_1,BIT10,BIT11,BIT12,BIT12_1,BFI14_0}. Such methods, based on the concept of {\em discrete line integral}  \cite{IP07, IP08, IT09}, have already been used in the context of Hamiltonian PDEs to derive the full discretization, when using the finite-difference MOL  approach in space \cite{BFI14}.

With this premise the paper is organized as follows: in Section~\ref{fourier} we study the case in which problem (\ref{wave}) is completed by the periodic boundary conditions (\ref{perbc}); in Section~\ref{sinexp} we study the case of general Dirichlet boundary condition (\ref{diribc}), whereas the case of general Neumann condition (\ref{neubc}) will be examined in Section~\ref{cosexp}; a few numerical tests are collected in Section~\ref{numtest} and, at last, some concluding remarks are given in Section~\ref{final}.

\section{The case of periodic boundary conditions}\label{fourier}
Let us consider the following complete set of orthonormal functions in $[0,1]$:
\begin{equation}\label{cnsn}
c_0(x)\equiv 1, \qquad c_k(x) = \sqrt{2}\cos(2k\pi x), \quad s_k(x) =\sqrt{2}\sin(2k\pi x), \qquad k=1,2,\dots,
\end{equation}
so that
\begin{equation}\label{orto}
\int_0^1 c_i(x)c_j(x)\dd x=\int_0^1 s_i(x)s_j(x)\dd x = \delta_{ij}, \qquad \int_0^1 c_i(x)s_j(x)\dd x=0, \qquad
\forall i,j,
\end{equation}
$\delta_{ij}$ being the Kronecker symbol. The following expansion of the solution of (\ref{wave})-(\ref{perbc}) is a slightly different way of writing the usual Fourier expansion in space:
\begin{eqnarray}\nonumber
u(x,t) &=& c_0(x)\gamma_0(t) +\sum_{n\ge 1} \left[c_n(x)\gamma_n(t)+s_n(x)\eta_n(t)\right] \\
\label{expu}
&\equiv&\gamma_0(t) +\sum_{n\ge 1}\left[ c_n(x)\gamma_n(t)+s_n(x)\eta_n(t)\right], \qquad x\in[0,1],\quad t\ge0,
\end{eqnarray}
with
$$
\gamma_n(t) = \int_0^1 c_n(x)u(x,t)\dd x, \qquad \eta_n(t) = \int_0^1 s_n(x)u(x,t)\dd x,
$$
which is allowed because of the periodic boundary conditions (\ref{perbc}). Consequently, by taking into account (\ref{orto}),  the first equation in (\ref{wave}) can be rewritten as: 
\begin{eqnarray}\nonumber
\ddgam_n(t) &=& -\aa^2(2\pi n)^2 \gamma_n(t) \\ \nonumber
&&- \int_0^1 c_n(x)f'\left( \gamma_0(t) +\sum_{n\ge 1} \left[c_n(x)\gamma_n(t)+s_n(x)\eta_n(t)\right] \right)\dd x, \quad n\ge0,\\ \label{fourier1} \\ \nonumber
\ddeta_n(t) &=& -\aa^2(2\pi n)^2 \eta_n(t)\\ \nonumber
&& - \int_0^1 s_n(x)f'\left( \gamma_0(t) +\sum_{n\ge 1}\left[ c_n(x)\gamma_n(t)+s_n(x)\eta_n(t)\right]  \right)\dd x, \quad n\ge1,
\end{eqnarray}
where the dot denotes, as usual, the time derivative. 
The initial conditions are clearly given by (see (\ref{wave})):
\begin{eqnarray}\nonumber
\gamma_n(0) = \int_0^1 c_n(x) \psi_0(x)\dd x, &\qquad& \eta_n(0) = \int_0^1 s_n(x) \psi_0(x)\dd x,\\[-1.5mm] \label{gameta0}\\[-1.5mm]
\dot\gamma_n(0) = \int_0^1 c_n(x) \psi_1(x)\dd x, &\qquad& \dot\eta_n(0) = \int_0^1 s_n(x) \psi_1(x)\dd x.\nonumber
\end{eqnarray}
By introducing the infinite vectors
\begin{eqnarray}\nonumber
\bom(x) &=& \pmatrix{cccccc} c_0(x)& c_1(x)& s_1(x) & c_2(x)& s_2(x)& \dots\endpmatrix^\top,\\
\label{qp} \\ \nonumber
\bq(t) &=& \pmatrix{cccccc} \gamma_0(t) & \gamma_1(t)& \eta_1(t) & \gamma_2(t)& \eta_2(t)& \dots\endpmatrix^\top,
\end{eqnarray}
the infinite matrix
\begin{equation}\label{D}
D = \pmatrix{cccccc} 0 \\ &(2\pi)^2 \\ &&(2\pi)^2\\ &&&(4\pi)^2\\ &&&&(4\pi)^2\\ &&&&& \ddots\endpmatrix,
\end{equation} and considering that (see (\ref{expu}))
\begin{equation}\label{expu1}
u(x,t) = \bom(x)^\top \bq(t),
\end{equation}
problem (\ref{fourier1}) can be cast in vector form as:
\begin{eqnarray} \label{fourier2} 
\dot\bq(t) &=& \bp(t), \qquad t>0, \\ \nonumber
\dot\bp(t) &=& -\aa^2D\bq(t) - \int_0^1 \bom(x) f'(\bom(x)^\top \bq(t))\dd x.
\end{eqnarray}
The following result holds true.

\begin{theo}\label{thH1} Problem (\ref{fourier2}) is Hamiltonian, with Hamiltonian
\begin{equation}\label{H2}
H(\bq,\bp) = \frac{1}2\bp^\top \bp + \frac{\aa^2}2\bq^\top D\bq + \int_0^1 f(\bom(x)^\top \bq)\dd x.
\end{equation}
This latter is equivalent to the Hamiltonian (\ref{v})-(\ref{H1}), via the expansion (\ref{expu})-(\ref{expu1}).
\end{theo}
\proof The first statement is straightforward, by considering that
$$\nabla_\bq f(\bom(x)^\top \bq)) = f'(\bom(x)^\top \bq)\bom(x).$$
The second statement then easily follows, by taking into account (\ref{expu1}), from the fact that, see (\ref{v}), (\ref{orto}), (\ref{expu}),  and (\ref{qp}):
\begin{eqnarray*}
\int_0^1 v(x,t)^2 \dd x &=& \int_0^1 u_t(x,t)^2\dd x ~=~ \int_0^1 \left(\dgam_0(t)+ \sum_{n\ge1}\left[ \dgam_n(t) c_n(x) + \deta_n(t) s_n(x)\right]\right)^2\dd x\\ &=& \dgam_0(t)^2+\sum_{n\ge1}\left[ \dgam_n(t)^2 +\deta_n(t)^2\right] ~\equiv\bp(t)^\top \bp(t),
\end{eqnarray*}
and
\begin{eqnarray*}
\int_0^1 u_x(x,t)^2 \dd x  &=& \int_0^1  \left(\sum_{n\ge1} 2\pi n\left[\eta_n(t) c_n(x) - \gamma_n(t) s_n(x)\right]\right)^2\dd x\\
&=& \sum_{n\ge1} (2\pi n)^2\left[\eta_n(t)^2+ \gamma_n(t)^2\right] ~=~ \bq(t)^\top D\bq(t).
\end{eqnarray*}
\QED

\subsection{Truncated Fourier approximation}\label{trunc}
In the computational practice, it is mandatory to truncate the infinite expansion (\ref{expu}) to a finite sum:\footnote{We refer, e.g., to \cite{CHQZ88}, for a corresponding comprehensive error analysis.}
\begin{equation}\label{expuN}
u(x,t) ~\approx~ \gamma_0(t) +\sum_{n=1}^N \left[c_n(x)\gamma_n(t)+s_n(x)\eta_n(t)\right] ~\equiv~ u_N(x,t).
\end{equation}
This reflects in the fact that the differential equations (\ref{fourier1}) now reduce to a finite number, i.e., $2N+1$. Correspondingly, one defines the finite vectors (see (\ref{qp})) in $\RR^{2N+1}$,
\begin{eqnarray}\nonumber
\bom_N(x) &=& \pmatrix{cccccccc} c_0(x)& c_1(x)& s_1(x) & c_2(x)& s_2(x)& \dots & c_N(x)& s_N(x)\endpmatrix^\top,\\ 
\label{qpN}\\[-3mm] \nonumber  
\bq_N(t) &=& \pmatrix{cccccccc} \gamma_0(t)& \gamma_1(t)& \eta_1(t) & \gamma_2(t)& \eta_2(t)& \dots & \gamma_N(t)& \eta_N(t)\endpmatrix^\top,\end{eqnarray}
and the matrix
\begin{equation}\label{DN}
D_N = \pmatrix{cccccccc} 0 \\ &(2\pi)^2 \\ &&(2\pi)^2\\ &&&(4\pi)^2\\ &&&&(4\pi)^2\\ &&&&&\ddots \\ &&&&&& (2N\pi)^2\\ &&&&&&& (2N\pi)^2\endpmatrix \in\RR^{2N+1\times 2N+1}.
\end{equation} Then, considering that (see (\ref{expuN}))
\begin{equation}\label{expu1N}
u_N(x,t) = \bom_N(x)^\top \bq_N(t),
\end{equation}
the equation which has to be satisfied by (\ref{expu1N}) can be cast in vector form as:
\begin{eqnarray} \label{fourier2N}
\dot\bq_N(t) &=& \bp_N(t), \qquad t>0, \\ \nonumber 
\dot\bp_N(t) &=& -\aa^2D_N\bq_N(t) - \int_0^1 \bom_N(x) f'(\bom_N(x)^\top \bq_N(t))\dd x,
\end{eqnarray}
for a total of $4N+2$ differential equations. Clearly, from (\ref{gameta0}) one obtains that the initial conditions for (\ref{fourier2N}) are given by:
$$
\bq_N(0) = \int_0^1 \bom_N(x)\psi_0(x)\dd x, \qquad  \bp_N(0) = \int_0^1 \bom_N(x)\psi_1(x)\dd x.
$$
The following result then easily follows by means of arguments similar to those used to prove Theorem~\ref{thH1}.

\begin{theo}\label{thH2} Problem (\ref{fourier2N}) is Hamiltonian, with Hamiltonian
\begin{equation}\label{H2N}
H_N(\bq_N,\bp_N) = \frac{1}2\bp_N^\top \bp_N + \frac{\aa^2}2\bq_N^\top D_N\bq_N + \int_0^1 f(\bom_N(x)^\top \bq_N)\dd x.
\end{equation}
This latter is equivalent to a truncated Fourier expansion of the Hamiltonian (\ref{v})-(\ref{H1}) (see also (\ref{H2})), that is, by truncating the expansion (\ref{expu})-(\ref{expu1}) as done in (\ref{expuN})-(\ref{expu1N}).
\end{theo}

\subsection{Approximating the integrals in space}\label{intx}
Clearly, the integral appearing in (\ref{fourier2N}) need to be, in general, approximated by means of a suitable quadrature rule. For this purpose, it could be convenient to do this by means of a composite trapezoidal rule, due to the fact that the argument is a periodic function. Consequently, having set
\begin{equation}\label{g_N}
g_N(x,t) = \bom_N(x) f'(\bom_N(x)^\top \bq_N(t)),
\end{equation} 
the uniform mesh on $[0,1]$,
\begin{equation}\label{xi}
x_i = i\Delta x, \quad i=0,\dots,m,\quad \Delta x = \frac{1}m,
\end{equation}
and considering that 
$$g_N(0,t) = g_N(1,t),$$
one obtains:
\begin{eqnarray}\nonumber
\int_0^1 g_N(x,t)\dd x &=& \Delta x\sum_{i=1}^m \frac{g_N(x_{i-1},t)+g_N(x_i,t)}2 ~+~R(m)\\ 
&=& \frac{1}m\sum_{i=0}^{m-1} g_N(x_i,t) ~+~ R(m).\label{compo}
\end{eqnarray}
Let us study the error $R(m)$. For this purpose, we need some preliminary result.

\begin{lem}\label{lem1}
Let us consider the trigonometric polynomial 
\begin{equation}\label{poly}
p(x) = \sum_{k=0}^K \left[a_k\cos(2k\pi x)+ b_k \sin(2k\pi x)\right],
\end{equation}
and the uniform mesh (\ref{xi}). Then, for all $m\ge K+1$, one obtains: $$\int_0^1 p(x)\dd x = \frac{1}m\sum_{i=0}^{m-1} p(x_i).$$
\end{lem}
\proof See, e.g., \cite[Th.\,5.1.4]{DaBi08}.\,\QED
\bigskip

\begin{lem}\label{lem2}
Let us consider the trigonometric polynomial (\ref{poly}) and the uniform mesh (\ref{xi}).
Then,  for all $m\ge N+K+1$, one obtains: 
\begin{eqnarray}\label{cosint}
\int_0^1 \cos(2j\pi x) p(x)\dd x &=& \frac{1}m\sum_{i=0}^{m-1} \cos(2j\pi x_i) p(x_i),\\ \label{sinint}
\int_0^1 \sin(2j\pi x) p(x)\dd x &=& \frac{1}m\sum_{i=0}^{m-1} \sin(2j\pi x_i) p(x_i),\qquad j=0,\dots,N.
\end{eqnarray}
\end{lem}
\proof
By virtue of the prosthaphaeresis formulae, one has, for all $j=0,\dots,N$ and $k=0,\dots,K$:
\begin{eqnarray*}
\cos(2j\pi x)\cos(2 k\pi x) &=& \frac{1}2\left[ \cos(2(k+j)\pi x) + \cos(2(k-j)\pi x)\right],\\
\cos(2j\pi x)\sin(2 k\pi x)  &=&  \frac{1}2\left[ \sin(2(k+j)\pi x) + \sin(2(k-j)\pi x) \right],\\
\sin(2j\pi x)\cos(2 k\pi x) &=& \frac{1}2\left[ \sin(2(k+j)\pi x) - \sin(2(k-j)\pi x)\right],\\
\sin(2j\pi x)\sin(2 k\pi x)  &=&  \frac{1}2\left[ \cos(2(k-j)\pi x) - \cos(2(k+j)\pi x) \right].
\end{eqnarray*}
Consequently, the integrals at the left-hand side in (\ref{cosint})--(\ref{sinint}) are trigonometric polynomials of degree at most $N+K$. By virtue of Lemma~\ref{lem1}, it then follows that they are exactly computed by means of the composite trapezoidal rule at the corresponding right-hand sides, provided that $m\ge N+K+1$.\QED 
\bigskip

By virtue of Lemma~\ref{lem2}, the following result follows at once.

\begin{theo}\label{esatti}
Let the function $f$ appearing in (\ref{g_N}) (see also (\ref{expu1N})) be a polynomial of degree $\nu$, and let us consider the uniform mesh (\ref{xi}). Then, with reference to (\ref{compo}), for all $m\ge \nu N+1$ one obtains:
$$R(m)= 0\qquad i.e.,\qquad \int_0^1 g_N(x,t)\dd x = \frac{1}m\sum_{i=0}^{m-1} g_N(x_i,t).$$
\end{theo}

For a general function $f$, the following result holds true.

\begin{theo}\label{approssimati}
Let the function $g_N(x,t)$ defined at (\ref{g_N}), with $t$ a fixed parameter, belong to $W_{per}^{r,p}$, the Banach space of periodic functions on $\RR$ whose distribution derivatives up to order $r$ belong to $L_{per}^p(\RR)$. Then, with reference to (\ref{xi})-(\ref{compo}), one has:
$$ R(m)  = O(m^{-r}).$$
\end{theo}
\proof See \cite[Th.\,1.1]{KuRa09}.\QED\bigskip

We end this section by mentioning that different approaches could be used, for approximating the integral appearing in (\ref{fourier2N}): as an example, we refer to \cite{EvWe99}, for a comprehensive review on this topic.

\subsection{Time integration}\label{HBVM}

Since problem (\ref{fourier2N}) is, for all $N\ge0$, Hamiltonian of dimension $4N+2$,  with Hamiltonian (\ref{H2N}), it is appropriate the use of an energy-conserving method for its numerical solution. We shall here consider, in particular, the family of Runge-Kutta type methods named {\em Hamiltonian Boundary Value Methods (HBVMs)} \cite{BIT09,BIT09_1,BIT10,BIT12,BIT12_1} (see also \cite{BIT11,BFI14_0}), already considered in \cite{BFI14}. Such methods rely on the concept of {\em discrete line integral}, introduced in \cite{IP07,IP08,IT09}, which is the discrete counterpart of the line integral for conservative vector fields. In particular, a HBVM$(k,s)$ method is the $k$-stages Runge-Kutta method, with $k\ge s$, defined by the following Butcher tableau:
\begin{equation}\label{RK}
\begin{array}{c|c}
\bc & A \equiv \P_{s+1}\hX_s\P_s^\top \Omega\\
\hline
       & \bb^\top 
       \end{array}
       \end{equation}
where the vectors $$\bc=(c_1,\dots,c_k)^\top ,    \qquad \bb = (b_1,\dots,b_k)^\top ,$$
contain the nodes and weights of the Gauss-Legendre formula of order $2k$, respectively, 
$$\Omega = \pmatrix{ccc} b_1\\ &\ddots \\ && b_k\endpmatrix,$$ $$\hX_s = \pmatrix{cccc}
\frac{1}2 &-\xi_1\\
\xi_1      &0         &\ddots\\
              &\ddots &\ddots &-\xi_{s-1}\\
              &            &\xi_{s-1} & 0\\
              \hline
              &            &               &\xi_s\endpmatrix \equiv \pmatrix{c} X_s \\
              \hline
              0\,\dots\,0~\xi_s\endpmatrix, \qquad  \xi_i = (4i^2-1)^{-\frac{1}2},
              $$
and, by setting
$\{P_j\}_{j\ge0}$ the family of Legendre polynomials, shifted and scaled so that
$$\int_0^1 P_i(c)P_j(c)\dd c = \delta_{ij}, \qquad \forall i,j=0,1,\dots,$$   
matrices $\P_s$ and $\P_{s+1}$ are defined as
$$\P_r = \left( P_{j-1}(c_i) \right) \in\RR^{k\times r}, \qquad r=s,s+1.$$
In particular, when $k=s$, (\ref{RK}) reduces to
$$
\begin{array}{c|c}
\bc & \P_sX_s\P_s^{-1}\\
\hline
       & \bb^\top 
       \end{array}$$
i.e., the Butcher tableau of the $s$-stage Gauss-Legendre collocation method. For this reason,
(\ref{RK}) can be also thought of as a generalization of the $W$-transform, as defined in \cite[page\,79]{HaWa96}.
The following result holds true \cite{BIT12_1}.

\begin{theo}\label{hbvmks}
For all $k\ge s$, the HBVM$(k,s)$ method (\ref{RK}), when applied for solving a Hamiltonian problem with stepsize $h$:
\begin{itemize}
\item is symmetric; \item has order $2s$;
\item is energy conserving when applied to polynomial Hamiltonians of degree
$\nu\le \lfloor \frac{2k}s\rfloor$;
\item for general and suitably regular Hamiltonians, the energy error at each step is $O(h^{2k+1})$.
\end{itemize}
\end{theo}

\begin{rem} From the result of the previous Theorem~\ref{hbvmks}, one has that an (at least {\em practical}) energy-conservation can be gained, for suitably regular Hamiltonians, provided that $k$ is large enough. On the other hand, this is not a big issue, from a computational point of view. In fact,  it turns out that the computational cost of a HBVM$(k,s)$ essentially depends on $s$. As a matter of fact, the discrete problem generated by the method can be seen to have dimension $s$, independently of $k$ \cite{BIT09,BIT12_1}. This fact, in turn, allows for an efficient implementation of the methods \cite{BIT11,BFI13,BFI14_0}.
\end{rem}

\section{The case of Dirichlet boundary conditions}\label{sinexp}
Let us now consider the case when the evolution equations are coupled with Dirichlet boundary conditions, so that the problem at hand is given by (\ref{wave})-(\ref{diribc}). There are several ways to cope with it: we shall sketch a couple of them in the subsections below.

\subsection{First approach}\label{firstapp}

A straightforward approach, quite easy to implement, is given by considering the auxiliary function  
\begin{equation}\label{zetag}
z(x,t) = u(x,t)-\left(x-\frac{1}2\right)\left[\phi_1(t) -\phi_0(t)\right] \equiv u(x,t) - g(x,t).
\end{equation}
In fact, the following result holds true.

\begin{theo}\label{perd}
Let $u(x,t)$ be the solution of problem (\ref{wave})-(\ref{diribc}). Then $z(x,t)$, defined at (\ref{zetag}), is the solution of the following problem with periodic boundary conditions.
\begin{eqnarray}\nonumber 
z_{tt}(x,t)&=&\aa^2z_{xx}(x,t)-f'(z(x,t)+g(x,t))-g_{tt}(x,t), \qquad (x,t)\in (0,1)\times(0,\infty),\\  \nonumber
z(x,0)&=& \psi_0(x) - g(x,0) ~\equiv~ \Psi_0(x), \\ \label{wavedsec}
z_t(x,0) &=& \psi_1(x) - g_t(x,0) ~\equiv~ \Psi_1(x), \qquad x\in(0,1),  \\ \nonumber
z(0,t) &=& z(1,t), \qquad\qquad t>0. 
\end{eqnarray}
\end{theo}
\proof
In fact, the first three equations in (\ref{wavedsec}) easily follow from (\ref{wave}) and (\ref{zetag}).
Moreover, because of the compatibility conditions,
$$\psi_0(0)=\phi_0(0), \quad \psi_0(1)=\phi_1(0), \quad \psi_1(0)=\phi_0'(0), \quad \psi_1(1)=\phi_1'(0),$$
one derives that
$$g(0,0) = -g(1,0) = \frac{\phi_0(0)-\phi_1(0)}2 \qquad \Rightarrow \qquad z(0,0) = z(1,0) = \frac{\phi_0(0)+\phi_1(0)}2,$$
$$g_t(0,0) = -g_t(1,0) = \frac{\phi_0'(0)-\phi_1(0)'}2 \qquad \Rightarrow \qquad z_t(0,0) = z_t(1,0) = \frac{\phi_0'(0)+\phi_1'(0)}2,$$
i.e., the initial conditions in (\ref{wavedsec}) are periodic. The thesis competes by observing that
$$g(0,t) = -g(1,t) = \frac{\phi_0(t)-\phi_1(t)}2 \qquad \Rightarrow \qquad z(0,t) = z(1,t) = \frac{\phi_0(t)+\phi_1(t)}2.\QED$$

\bigskip Based on this result, by using the same notation as in (\ref{qp})--(\ref{D}), we can then look for a Fourier expansion in the form (compare with (\ref{expu1}))
\begin{equation}\label{expuz}
z(x,t) = \bom(x)^\top\bq(t),
\end{equation}
thus arriving at the infinite set of differential equations
\begin{eqnarray}\label{fourierz}
\dot\bq(t) &=& \bp(t),  \qquad t>0,\\  \nonumber
\dot\bp(t) &=& -\aa^2D\bq(t) - \int_0^1 \bom(x) \left[ f'(\bom(x)^\top \bq(t) +g(x,t)) +g_{tt}(x,t)\right]\dd x.
\end{eqnarray}
The following result holds true, whose proof is similar to that of Theorem~\ref{thH1}.

\begin{theo}
Problem (\ref{fourierz}) is Hamiltonian, with  non-autonomous Hamiltonian
\begin{eqnarray*}\nonumber
\lefteqn{H(\bq,\bp,t) ~=~ \frac{1}2\bp^\top \bp + \frac{\aa^2}2\bq^\top D\bq }\\ \label{Hz}
&&+ \int_0^1 \left[f(\bom(x)^\top \bq +g(x,t))+g_{tt}(x,t)\, \bom(x)^\top \bq\right]\dd x.
\end{eqnarray*}
\end{theo}

A finite-dimensional approximation of (\ref{expuz}) can then be derived by using similar arguments as those seen in Section~\ref{trunc}. In more details, by using the notation (\ref{qpN})-(\ref{DN}), one looks for a truncated Fourier expansion in the form (compare with (\ref{expu1N})):
$$
z_N(x,t) = \bom_N(x)^\top\bq_N(t),
$$
thus arriving at the following set of $2(2N+1)$ differential equations:
\begin{eqnarray}\label{fourierzN} 
\dot\bq_N(t) &=& \bp_N(t),  \qquad t>0,\\ \nonumber
\dot\bp_N(t) &=& -\aa^2D_N\bq_N(t) - \int_0^1 \bom_N(x) \left[ f'(\bom_N(x)^\top \bq_N(t) +g(x,t)) +g_{tt}(x,t)\right]\dd x.
\end{eqnarray}
The following result then easily follows.

\begin{theo}
Problem (\ref{fourierzN}) is Hamiltonian, with non-autonomous Hamiltonian
\begin{eqnarray*}\nonumber
\lefteqn{H_N(\bq,\bp,t)~=~ \frac{1}2\bp_N^\top \bp_N + \frac{\aa^2}2\bq_N^\top D_N\bq_N }\\ \label{HzN} 
&& + \int_0^1 \left[f(\bom_N(x)^\top \bq_N +g(x,t)) +g_{tt}(x,t)\, \bom_N(x)^\top \bq_N\right]\dd x.
\end{eqnarray*}
Moreover, along the solution of (\ref{fourierzN}),
\begin{eqnarray}\nonumber
\lefteqn{\dot{H}_N(\bq,\bp,t) ~\equiv~\frac{\partial}{\partial t}H_N(\bq,\bp,t) }\\ \label{derHN}
&=&\int_0^1 \left[f'(\bom_N(x)^\top \bq_N +g(x,t))g_t(x,t) +g_{ttt}(x,t)\, \bom_N(x)^\top \bq_N\right]\dd x.
\end{eqnarray}
\end{theo}

\begin{rem}
The main difference, with respect to the case of periodic boundary conditions studied in Section~\ref{trunc}, stems from the fact that now the Hamiltonian is time dependent. Moreover, one has to consider that the involved integrals have to be, in general, approximated by means of different quadrature rules (e.g., a high-order composite Newton-Cotes or Gaussian formula), than those exposed in Section~\ref{intx}, due to the fact that now, in general, the integrand is no more a periodic function in the space argument. We omit, however, the details about this standard argument.
\end{rem}

It is worth noting that, following the approach in \cite{BFI14}, we can ``embed'' problem (\ref{fourierzN}) (as well as its infinite counterpart (\ref{fourierz})) into a higher dimensional Hamiltonian problem, with an autonomous Hamiltonian. In fact, by introducing the auxiliary scalar conjugate variables $\hq$ and $\hp$, and the augmented (autonomous) Hamiltonian
\begin{eqnarray}\nonumber
\lefteqn{ \hat{H}_N(\bq_N,\bp_N,\hq,\hp) ~=~ \frac{1}2\bp_N^\top \bp_N + \frac{\aa^2}2\bq_N^\top D_N\bq_N}  \\
&&+ \int_0^1 \left[f(\bom_N(x)^\top \bq_N +g(x,\hq)) +g_{\hq\hq}(x,\hq)\, \bom_N(x)^\top \bq_N\right]\dd x 
+\hp \nonumber \\
&\equiv& H_N(\bq_N,\bp_N,\hq) + \hp, \label{augHz}
\end{eqnarray}
one obtains the augmented Hamiltonian problem (see (\ref{derHN}))\begin{eqnarray}\nonumber
\dot\bq_N(t) &=& \bp_N(t), \\ \nonumber
\dot\bp_N(t) &=& -\aa^2D_N\bq_N(t) - \int_0^1 \bom_N(x) \left[ f'(\bom_N(x)^\top \bq_N(\hq) +g(x,t)) +g_{\hq\hq}(x,\hq)\right]\dd x,\\
\dot\hq &=& 1,    \label{fourierzN1} \\ \nonumber
\dot\hp &=& -\frac{\partial}{\partial \hq} H_N(\bq_N,\bp_N,\hq), \qquad\qquad t>0.
\end{eqnarray}
By using the initial conditions (see (\ref{wavedsec}))\,\footnote{As is clear, in (\ref{fourierzN1})-(\ref{qp0}) $\hq(t)\equiv t$.}
\begin{equation}\label{qp0}
\bq_N(0) = \int_0^1\bom_N(x)\Psi_0(x)\dd x, \quad \bp_N(0) = \int_0^1\bom_N(x)\Psi_1(x)\dd x,\quad \hq(0)=\hp(0)=0, 
\end{equation}
the following straightforward result easily follows (see, e.g., \cite{BFI14}).
\begin{theo}
Along the solution of (\ref{fourierzN1})-(\ref{qp0}), one has
\begin{equation}\label{constaugH}
\hat{H}_N(\bq_N(t),\bp_N(t),\hq(t),\hp(t)) \equiv \hat{H}_N(\bq_N(0),\bp_N(0),0,0) \equiv H_N(\bq_N(0),\bp_N(0),0),
\end{equation}
for all ~$t\ge0$.
\end{theo}
\begin{rem} Clearly, a suitable HBVM$(k,s)$ formula can be conveniently used for numerically solving (\ref{fourierzN1})-(\ref{qp0}), and fulfilling, at least ``practically'', (\ref{constaugH}), accordingly with the results of Theorem~\ref{hbvmks}.\end{rem}

\subsection{A second approach}\label{secondapp}

Another approach for solving (\ref{wave})-(\ref{diribc}) is obtained by considering the following associated linear problem,
\begin{eqnarray}\nonumber 
\hu_{tt}(x,t)&=&\aa^2\hu_{xx}(x,t), \qquad (x,t)\in (0,1)\times(0,\infty),\\  \nonumber
\hu(x,0)&=&\psi_0(x),\\  \label{waveldiri}
\hu_t(x,0)&=&\psi_1(x), \qquad x\in(0,1), \\ \nonumber
\hu(0,t) &=& \phi_0(t), \\ \nonumber
\hu(1,t) &=& \phi_1(t), \qquad t>0,
\end{eqnarray}
whose solution we assume to know (a detailed discussion is presented in Section~\ref{aux1} below). Let us then define the auxiliary function
\begin{equation}\label{zeta}
z(x,t) = u(x,t)-\hu(x,t).
\end{equation}
It is straightworfard to check that it satisfies the non-autonomous nonlinear wave problem: 
\begin{eqnarray}\nonumber 
z_{tt}(x,t)&=&\aa^2z_{xx}(x,t)-f'(z(x,t)+\hu(x,t)), \qquad (x,t)\in (0,1)\times(0,\infty),\\  \nonumber
z(x,0)&=& z_t(x,0) ~=~0, \qquad x\in(0,1), \\ \label{waved0}
z(0,t) &=& z(1,t) ~=~0, \qquad t>0,
\end{eqnarray}
whose solution put in the form (compare with (\ref{cnsn})-(\ref{orto}))
\begin{equation}\label{zsin}
z(x,t) = \sum_{n\ge1} \hs_n(x)\eta_n(t), \qquad x\in[0,1], \quad t\ge0,
\end{equation}
where we are now considering the orthonormal basis, on [0,1], of the continuous functions which vanish at the end-points of the interval,
\begin{equation}\label{soloseni}
\hs_i(x) = \sqrt{2}\sin(i\pi x), \qquad \int_0^1\hs_i(x)\hs_j(x)\dd x = \delta_{ij}, \qquad \forall i,j\ge1.
\end{equation}
Consequently, (\ref{zsin}) satisfies the homogeneous boundary conditions in (\ref{waved0}). Moreover, because of the initial conditions in (\ref{waved0}), one obtains the (infinite) differential problem:
\begin{eqnarray}\nonumber
\ddeta_n(t) &=& -\aa^2(\pi n)^2 \eta_n(t) - \int_0^1\hs_n(x)f'( z(x,t)+\hu(x,t))\dd x, \qquad t>0,\\
\eta_n(0) &=& \deta_n(0)~=~0, \qquad n\ge1.\label{eta0} 
\end{eqnarray}
As done in the case of periodic boundary conditions, we can cast this problem in vector form by defining the (infinite) vectors (compare with (\ref{qp}))\,\footnote{In order to emphasize the similarities, also avoiding to introduce a more involved notation, we shall use the same notation used in Section~\ref{firstapp}.}

\begin{equation}\label{qp1}
\bom(x) = \pmatrix{ccc} \hs_1(x)& \hs_2(x)& \dots\endpmatrix^\top, \qquad
\bq(t) = \pmatrix{ccc} \eta_1(t)& \eta_2(t)& \dots\endpmatrix^\top, 
\end{equation} and the the infinite matrix (compare with (\ref{D}))
\begin{equation}\label{D1}
D = \pmatrix{ccc} \pi^2\\ & (2\pi)^2\\ &&\ddots\endpmatrix,
\end{equation}
so that (compare with (\ref{expu1}))
\begin{equation}\label{expz}
z(x,t) = \bom(x)^\top \bq(t).
\end{equation}
Consequently, (\ref{eta0}) can be cast in Hamiltonian form, by taking into account (\ref{expz}) (compare with (\ref{fourier2})), as
\begin{eqnarray}\nonumber 
\dbq(t) &=& \bp(t),  \qquad t>0, \\ \label{fourierD}
\dbp(t) &=& -\aa^2D\bq(t) -\int_0^1\bom(x)f'(\bom(x)^\top \bq(t)+\hu(x,t))\dd x,\\ \nonumber
\bq(0) &=& \bp(0) ~=~ \bf0,
\end{eqnarray}
with non-autonomous Hamiltonian
\begin{equation}\label{Hpqt}
H(\bq,\bp,t) = \frac{1}2\bp^\top \bp +\frac{\aa^2}2\bq D\bq +\int_0^1f(\bom(x)^\top \bq+\hu(x,t))\dd x.
\end{equation}
In a similar way as it has been done in Section~\ref{trunc} for the case of periodic boundary conditions, one derives a practical procedure by approximating the infinite expansion (\ref{zsin}) through a truncated one,
\begin{equation}\label{expzN}
z(x,t) \approx z_N(x,t) = \sum_{n=1}^N \hs_n(x)\eta_n(t) \equiv \bom_N(x)^\top \bq_N(t),
\end{equation}
with (compare with (\ref{qpN})-(\ref{DN}))
\begin{equation}\label{qpzN}
\bom_N(x) = \pmatrix{c} \hs_1(x)\\ \vdots \\ \hs_N(x)\endpmatrix, \quad
\bq_N(t) = \pmatrix{c} \eta_1(t)\\ \vdots \\ \eta_N(t)\endpmatrix, \quad
D_N =  \pmatrix{ccc} \pi^2\\ &\ddots \\ && (N\pi)^2\endpmatrix,
\end{equation}
so that we arrive at the Hamiltonian problem (compare with (\ref{fourier2N}))
\begin{eqnarray}\nonumber
\dbq_N(t) &=& \bp_N(t), \qquad \qquad t>0,\\ \nonumber
\dbp_N(t) &=& -\aa^2D_N\bq_N(t) -\int_0^1\bom_N(x)f'(\bom_N(x)^\top \bq_N(t)+\hu(x,t))\dd x, \\ \label{fourierN} 
\bq_N(0) &=& \bp_N(0) ~=~\bf0,
\end{eqnarray}
with  non-autonomous Hamiltonian (compare with (\ref{H2N}))
\begin{equation}\label{HN}
H_N(\bq_N,\bp_N,t) = \frac{1}2\bp_N^\top \bp_N +\frac{\aa^2}2\bq_N D_N\bq_N +\int_0^1f(\bom_N(x)^\top \bq_N+\hu(x,t))\dd x.
\end{equation}
Similarly as done in Section~\ref{firstapp}, problem (\ref{fourierN}) can be ``embedded'' in the higher-dimensional problem defined by the augmented (autonomous) Hamiltonian (compare with (\ref{augHz}))
\begin{equation}\label{augHz1}
\hat{H}_N(\bq_N,\bp_N,\hq,\hp) = H_N(\bq_N,\bp_N,\hq)+\hp,
\end{equation}
obtained by introducing the auxiliary scalar conjugate variables $\hq$ and $\hp$. This latter problem can then be conveniently solved by using a suitable HBVM$(k,s)$ formula.

\subsubsection{Solving the auxiliary linear problem}\label{aux1}
For solving the auxiliary linear problem (\ref{waveldiri}), we shall further consider the function
$$
\bu(x,t) = \hu(x,t) -x\phi_1(t)-(1-x)\phi_0(t),
$$
satisfying the following additional problem,
\begin{eqnarray}\nonumber 
\bu_{tt}(x,t)&=&\aa^2\bu_{xx}(x,t) -x\phi_1''(t)-(1-x)\phi_0''(t)\\ \nonumber
&\equiv& \aa^2\bu_{xx}(x,t)+\baf(x,t), \qquad (x,t)\in (0,1)\times(0,\infty),\\  \label{waveldiri1}
\bu(x,0)&=&\psi_0(x)-x\psi_0(1)-(1-x)\psi_0(0) ~\equiv~\bpsi_0(x),\\  \nonumber
\bu_t(x,0)&=&\psi_1(x)-x\psi_1(1)-(1-x)\psi_1(0) ~\equiv~\bpsi_1(x), \qquad x\in(0,1), \\ \nonumber
\bu(0,t) &=& \bu(1,t) ~=~0, \qquad t>0,
\end{eqnarray}
whose solution is easily seen to be obtained as superposition of the solutions of the following two problems:
\begin{eqnarray}\nonumber 
\bau_{tt}(x,t)&=&\aa^2\bau_{xx}(x,t), \qquad (x,t)\in (0,1)\times(0,\infty),\\  \label{waveldiri2}
\bau(x,0)&=&\bpsi_0(x),\\  \nonumber
\bau_t(x,0)&=&\bpsi_1(x), \qquad x\in(0,1), \\ \nonumber
\bau(0,t) &=& \bau(1,t) ~=~0, \qquad t>0,
\end{eqnarray}
and
\begin{eqnarray}\nonumber 
\bbu_{tt}(x,t)&=&\aa^2\bbu_{xx}(x,t) +\baf(x,t), \qquad (x,t)\in (0,1)\times(0,\infty),\\  \label{waveldiri3}
\bbu(x,0)&=&\bbu_t(x,0)~=~0, \qquad x\in(0,1), \\ \nonumber
\bbu(0,t) &=& \bbu(1,t) ~=~0, \qquad t>0,
\end{eqnarray}
i.e., 
\begin{equation}\label{bu12}
\bu(x,t) = \bau(x,t)+\bbu(x,t).
\end{equation}
The following result is then easily derived.
\begin{theo}
The solutions of problems (\ref{waveldiri2}) and (\ref{waveldiri3}) are given by:
\begin{eqnarray}\label{bu1}
\bau(x,t) &=& \frac{1}2\left( G(x-\aa t)+G(x+\aa t) + \aa^{-1}\int_{x-\aa t}^{x+\aa t} H(s)\dd s\right),\\
\bbu(x,t) &=& \frac{1}{2\aa}\int_0^t \int_{x-\aa(t-\tau)}^{x+\aa(t-\tau)} F(s,\tau)\dd s\dd\tau,\label{bu2}
\end{eqnarray}
where $F(x,t)$, $G(x)$, and $H(x)$ are the following periodic functions (see (\ref{waveldiri1})):
\begin{eqnarray}\nonumber
F(x,t) &=&\left\{\begin{array}{ccl} \baf(x,t), &\quad &\mbox{if}\quad x\in(0,1),\\
-\baf(-x,t), &&\mbox{if}\quad x\in(-1,0),\\
F(y,t), &&\mbox{if}\quad x=y+2n, \quad y\in(-1,1),\quad n\in\ZZ_0,\end{array}\right.\\ \label{oddr}
G(x) &=&\left\{\begin{array}{ccl} \bpsi_0(x), &\quad &\mbox{if}\quad x\in(0,1),\\
-\bpsi_0(-x), &&\mbox{if}\quad x\in(-1,0),\\
G(y), &&\mbox{if}\quad x=y+2n, \quad y\in(-1,1),\quad n\in\ZZ_0,\end{array}\right.\\ \nonumber
H(x) &=&\left\{\begin{array}{ccl} \bpsi_1(x), &\quad &\mbox{if}\quad x\in(0,1),\\
-\bpsi_1(-x), &&\mbox{if}\quad x\in(-1,0),\\
H(y), &&\mbox{if}\quad x=y+2n, \quad y\in(-1,1),\quad n\in\ZZ_0.\end{array}\right.
\end{eqnarray}
\end{theo}
Consequently, (\ref{bu12}) can be obtained by using suitable quadrature rules.

\begin{rem}\label{numer1} We also mention that problem (\ref{waveldiri}) could be solved numerically by using a high-order boundary value method (BVM, see \cite{BTbook} for more details on BVMs). In particular, the high order formulae in \cite{AmSg05} are particularly suited. In so doing, at each integration step, one first solves (\ref{waveldiri}) on the time window $[0,h]$, then proceeds with  the solution of (\ref{waved0}) on the same time domain. This basic procedure is thus iterated to cover the desired time interval.
\end{rem}

\section{The case of Neumann boundary conditions}\label{cosexp}
Let us now consider the case when Neumann boundary conditions are coupled with the equations, so that, the problem at hand is given by (\ref{wave})-(\ref{neubc}). Also in this case, one can use different approaches to cope with this problem. In the following sections, we sketch two of them, which are similar to those examined in Sections~\ref{firstapp} and \ref{secondapp}, respectively, for the case of Dirichlet boundary conditions. Deliberately, we shall use very similar notations as those used there, in order to emphasize the existing similarities.

\subsection{First approach}\label{firstapp1}
Let us consider the auxiliary function
\begin{equation}\label{zetag1}
z(x,t) = u(x,t) - \left[\frac{x^2}2[ \phi_1(t)-\phi_0(t)] +x\phi_0(t) \right]
\equiv u(x,t)-g(x,t), 
\end{equation}
for which the following straightforward result holds true.
\begin{theo}\label{perdx}
Let $u(x,t)$ be the solution of problem (\ref{wave})-(\ref{neubc}). Then $z(x,t)$, defined at (\ref{zetag1}), is the solution of the following problem with periodic boundary conditions.
\begin{eqnarray}\nonumber 
z_{tt}(x,t)&=&\aa^2z_{xx}(x,t)-f'(z(x,t)+g(x,t))-\left[g_{tt}(x,t)-\aa^2g_{xx}(x,t)\right], \\ \nonumber
 &\equiv& \aa^2z_{xx}(x,t)-f'(z(x,t)+g(x,t))-\Phi(x,t)\qquad (x,t)\in (0,1)\times(0,\infty),\\  \label{wavedsec1}
z(x,0)&=& \psi_0(x) - g(x,0) ~\equiv~ \Psi_0(x), \\ \nonumber
z_t(x,0) &=& \psi_1(x) - g_t(x,0) ~\equiv~ \Psi_1(x), \qquad x\in(0,1),  \\ \nonumber
z_x(0,t) &=& z_x(1,t) ~=~0, \qquad\qquad t>0. 
\end{eqnarray}
\end{theo}
Consequently, we can look for an expansion of $z_x(x,t)$ in the form (see (\ref{soloseni}))
\begin{equation}\label{zxsin}
z_x(x,t) = -\pi\sum_{n\ge1} n\hs_n(x)\gamma_n(t), \qquad x\in[0,1], \quad t\ge0,
\end{equation}
i.e.,
\begin{equation}\label{zcos}
z(x,t) = \sum_{n\ge0} \hc_n(x)\gamma_n(t), \qquad x\in[0,1], \quad t\ge0,
\end{equation}
where we are considering the orthonormal functions on [0,1]:
\begin{equation}\label{solocoseni}
\hc_0(x)\equiv1, \qquad \hc_i(x) = \sqrt{2}\cos(i\pi x), \qquad \int_0^1\hc_i(x)\hc_j(x)\dd x = \delta_{ij}, \qquad \forall i,j.
\end{equation}
One then obtains the (infinite) differential problem (see (\ref{wavedsec1})):
\begin{eqnarray*}
\ddgam_n(t) &=& -\aa^2(\pi n)^2 \gamma_n(t) 
 -\int_0^1\hc_n(x)\left[f'( z(x,t)+g(x,t))+\Phi(x,t)\right]\dd x, \qquad t>0,\\
\gamma_n(0) &=& \int_0^1\hc_n(x)\Psi_0(x)\dd x, \qquad \dgam_n(0) ~=~ \int_0^1\hc_n(x)\Psi_1(x)\dd x, \qquad n\ge0,
\end{eqnarray*}
which we can cast in vector form as 
\begin{eqnarray}\nonumber 
\dbq(t) &=& \bp(t),  \qquad t>0, \\ \label{fourierD1}
\dbp(t) &=& -\aa^2D\bq(t) -\int_0^1\bom(x)\left[f'(\bom(x)^\top \bq(t)+g(x,t))+\Phi(x,t)\right]\dd x,\\ \nonumber
\bq(0) &=& \int_0^1\bom(x)\Psi_0(x)\dd x, \qquad \bp(0) ~=~ \int_0^1\bom(x)\Psi_1(x)\dd x,
\end{eqnarray}
with non-autonomous Hamiltonian
$$
H(\bq,\bp,t) = \frac{1}2\bp^\top \bp +\frac{\aa^2}2\bq D\bq +\int_0^1\left[f(\bom(x)^\top \bq+g(x,t))+\Phi(x,t)\,\bom(x)^\top \bq\right]\dd x,
$$
where (see (\ref{solocoseni}))
\begin{equation}\label{qpD2}
\bom(x) = \pmatrix{c} \hc_0(x)\\ \hc_1(x)\\ \hc_2(x)\\ \vdots\endpmatrix, \qquad
\bq(t) = \pmatrix{c} \gamma_0(t)\\ \gamma_1(t)\\ \gamma_2(t)\\ \vdots\endpmatrix, \qquad
D = \pmatrix{cccc} 0 \\ &\pi^2\\ && (2\pi)^2\\ &&&\ddots\endpmatrix.
\end{equation}
Similarly as done before, a finite dimensional approximation is now obtained by considering
\begin{equation}\label{expzN1}
z(x,t) \approx z_N(x,t) = \sum_{n=0}^N \hc_n(x)\gamma_n(t) \equiv \bom_N(x)^\top \bq_N(t),
\end{equation}
and
\begin{equation}\label{qpzN1}
\bom_N(x) = \pmatrix{c} \hc_0(x) \\ \hc_1(x)\\ \vdots \\ \hc_N(x)\endpmatrix, \quad
\bq_N(t) = \pmatrix{c} \gamma_0(t) \\ \gamma_1(t)\\ \vdots \\ \gamma_N(t)\endpmatrix, \quad
D_N =  \pmatrix{cccc} 0 \\ &\pi^2\\ &&\ddots \\ &&& (N\pi)^2\endpmatrix,
\end{equation}
thus obtaining the Hamiltonian problem, of dimension $2(N+1)$,
\begin{eqnarray}  \label{fourierN1} 
\dbq_N(t) &=& \bp_N(t), \qquad \qquad t>0,\\ \nonumber
\dbp_N(t) &=& -\aa^2D_N\bq_N(t) -\int_0^1\bom_N(x)\left[f'(\bom_N(x)^\top \bq_N(t)+g(x,t))+\Phi(x,t)\right]\dd x, \\ \nonumber
\bq_N(0) &=& \int_0^1\bom_N(x)\Psi_0(x)\dd x, \qquad \bp_N(0) ~=~ \int_0^1\bom_N(x)\Psi_1(x)\dd x,
\end{eqnarray}
with  non-autonomous Hamiltonian
\begin{eqnarray}\label{HN1}
H_N(\bq_N,\bp_N,t) &=& \frac{1}2\bp_N^\top \bp_N +\frac{\aa^2}2\bq_N D_N\bq_N \\&&+\int_0^1\left[f(\bom_N(x)^\top \bq_N+g(x,t))+\Phi(x,t)\,\bom_N(x)^\top \bq_N\right]\dd x.\nonumber
\end{eqnarray}
Similarly as done in Section~\ref{firstapp}, problem (\ref{fourierN1}) can be ``embedded'' in the higher-dimen\-sional problem defined by the augmented (autonomous) Hamiltonian, obtained by introducing the auxiliary scalar conjugate variables $\hq$ and $\hp$, which is formally still given by (\ref{augHz1}), with $H_N$ given by (\ref{HN1}). This latter problem can then be conveniently solved by using a suitable HBVM$(k,s)$ formula, by considering that, also in the present case, the arguments in Section~\ref{intx} need to be suitably modified, for the approximation of the integrals in space, due to the fact that now $\Phi(x,t)$ cannot be assumed to be periodic in space. 

\subsection{A second approach}\label{secondapp1} 
We repeat here similar steps as those in Section~\ref{secondapp}, by considering the associated linear problem,
\begin{eqnarray}\nonumber 
\hu_{tt}(x,t)&=&\aa^2\hu_{xx}(x,t), \qquad (x,t)\in (0,1)\times(0,\infty),\\  \nonumber
\hu(x,0)&=&\psi_0(x),\\  \label{wavelneu}
\hu_t(x,0)&=&\psi_1(x), \qquad x\in(0,1), \\ \nonumber
\hu_x(0,t) &=& \phi_0(t), \\ \nonumber
\hu_x(1,t) &=& \phi_1(t), \qquad t>0,
\end{eqnarray}
whose solution we assume to know (see Section~\ref{aux2} below). Let us then define the auxiliary function, formally still given by (\ref{zeta}),
which satisfies the non-autonomous nonlinear wave problem: 
\begin{eqnarray*}\nonumber 
z_{tt}(x,t)&=&\aa^2z_{xx}(x,t)-f'(z(x,t)+\hu(x,t)), \qquad (x,t)\in (0,1)\times(0,\infty),\\  \nonumber
z(x,0)&=& z_t(x,0) ~=~0, \qquad x\in(0,1), \\ \label{waven0}
z_x(0,t) &=& z_x(1,t) ~=~0, \qquad t>0.
\end{eqnarray*}
Since $z_x(x,t)$ vanishes at the end-points, we can look for an expansion identical to (\ref{zxsin}), thus arriving at the same expansion (\ref{zcos}) for $z(x,t)$. Consequently, (\ref{expz})--(\ref{Hpqt}) continue formally to hold, by replacing (\ref{qp1})-(\ref{D1}) with (\ref{qpD2}), as well as (\ref{expzN})--(\ref{augHz1}), by replacing (\ref{qpzN}) with (\ref{qpzN1}).

\subsection{Solving the auxiliary linear problem}\label{aux2}
In order for solving the auxiliary linear problem (\ref{wavelneu}), we follow a procedure formally very similar to that studied in Section~\ref{aux1} for the Dirichlet case. Let us then consider the function
\begin{equation}\label{bun}
\bu(x,t) = \hu(x,t) -\frac{x}2\left( (2-x) \phi_0(t)+x\phi_1(t)\right),
\end{equation}
satisfying the following additional problem,
\begin{eqnarray}\nonumber 
\bu_{tt}(x,t)&=&\aa^2\left[\bu_{xx}(x,t) -\phi_0(t)+\phi_1(t)\right] -\frac{x}2\left((2-x)\phi_0''(t)+x\phi_1''(t)\right)\\ \nonumber
&\equiv& \aa^2\bu_{xx}(x,t)+\baf(x,t), \qquad (x,t)\in (0,1)\times(0,\infty),\\  \label{wavelneu1}
\bu(x,0)&=&\psi_0(x)-\frac{x}2\left((2-x)\psi_0'(0)+x\psi_0'(1)\right) ~\equiv~\bpsi_0(x),\\  \nonumber
\bu_t(x,0)&=&\psi_1(x)-\frac{x}2\left((2-x)\psi_1'(0)+x\psi_1'(1)\right) ~\equiv~\bpsi_1(x), \qquad x\in(0,1), \\ \nonumber
\bu_x(0,t) &=& \bu_x(1,t) ~=~0, \qquad t>0,
\end{eqnarray}
whose solution is easily seen to be obtained as superposition of the solutions of the following two problems:
\begin{eqnarray}\nonumber 
\bau_{tt}(x,t)&=&\aa^2\bau_{xx}(x,t), \qquad (x,t)\in (0,1)\times(0,\infty),\\  \label{wavelneu2}
\bau(x,0)&=&\bpsi_0(x),\\  \nonumber
\bau_t(x,0)&=&\bpsi_1(x), \qquad x\in(0,1), \\ \nonumber
\bau_x(0,t) &=& \bau_x(1,t) ~=~0, \qquad t>0,
\end{eqnarray}
and
\begin{eqnarray}\nonumber 
\bbu_{tt}(x,t)&=&\aa^2\bbu_{xx}(x,t) +\baf(x,t), \qquad (x,t)\in (0,1)\times(0,\infty),\\  \label{wavelneu3}
\bbu(x,0)&=&\bbu_t(x,0)~=~0, \qquad x\in(0,1), \\ \nonumber
\bbu_x(0,t) &=& \bbu_x(1,t) ~=~0, \qquad t>0.
\end{eqnarray}
I.e., (\ref{bu12}) is still formally valid.
The following result is then easily derived.
\begin{theo}
The solutions of problems (\ref{wavelneu2}) and (\ref{wavelneu3}) are formally still given by
(\ref{bu1}) and (\ref{bu2}), respectively,
where $F(x,t)$, $G(x)$, and $H(x)$ are the following periodic functions (see (\ref{wavelneu1})):\footnote{We observe that the functions in (\ref{evenr}) are obtained through {\em even reflection} of the original functions, whereas the functions in (\ref{oddr}) are obtained by {\em odd reflection}.}
\begin{eqnarray}\nonumber
F(x,t) &=&\left\{\begin{array}{ccl} \baf(x,t), &\quad &\mbox{if}\quad x\in(0,1),\\
\baf(-x,t), &&\mbox{if}\quad x\in(-1,0),\\
F(y,t), &&\mbox{if}\quad x=y+2n, \quad y\in(-1,1),\quad n\in\ZZ_0,\end{array}\right.\\ \label{evenr}
G(x) &=&\left\{\begin{array}{ccl} \bpsi_0(x), &\quad &\mbox{if}\quad x\in(0,1),\\
\bpsi_0(-x), &&\mbox{if}\quad x\in(-1,0),\\
G(y), &&\mbox{if}\quad x=y+2n, \quad y\in(-1,1),\quad n\in\ZZ_0,\end{array}\right.\\ \nonumber
H(x) &=&\left\{\begin{array}{ccl} \bpsi_1(x), &\quad &\mbox{if}\quad x\in(0,1),\\
\bpsi_1(-x), &&\mbox{if}\quad x\in(-1,0),\\  
H(y), &&\mbox{if}\quad x=y+2n, \quad y\in(-1,1),\quad n\in\ZZ_0.\end{array}\right.
\end{eqnarray}
\end{theo}
Consequently, also in such a case, the function $\bu(x,t)$ in (\ref{bun}), formally still given by (\ref{bu12}), can be obtained by using suitable quadrature rules.

Finally, we observe that the same arguments in Remark~\ref{numer1} apply to this case.

\section{Numerical tests}\label{numtest}

\begin{figure}[ph]
\centerline{\includegraphics[width=9cm,height=6cm]{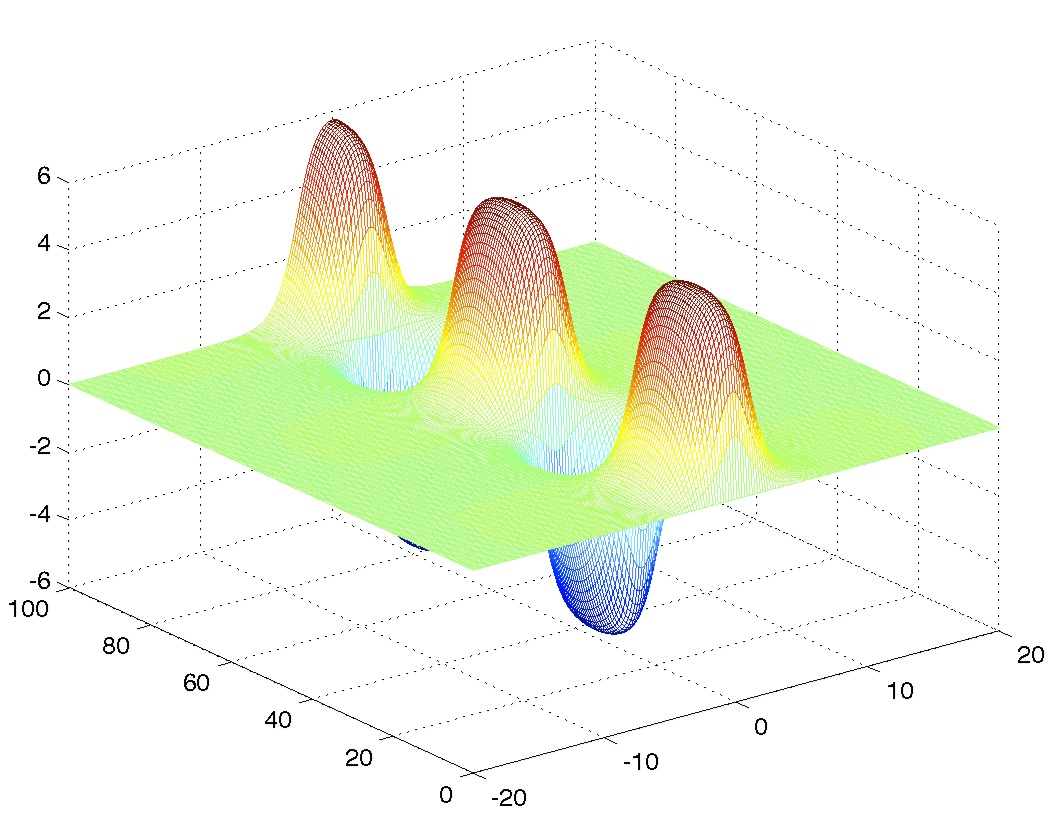}}
\caption{{\em Breather}, i.e.,  soliton-like solution (\ref{sineGu})-(\ref{fi}) of problem (\ref{sineG})-(\ref{sineG0}), $\gamma=1.01$.}
\label{soli0}

\medskip
\centerline{\includegraphics[width=9cm,height=6cm]{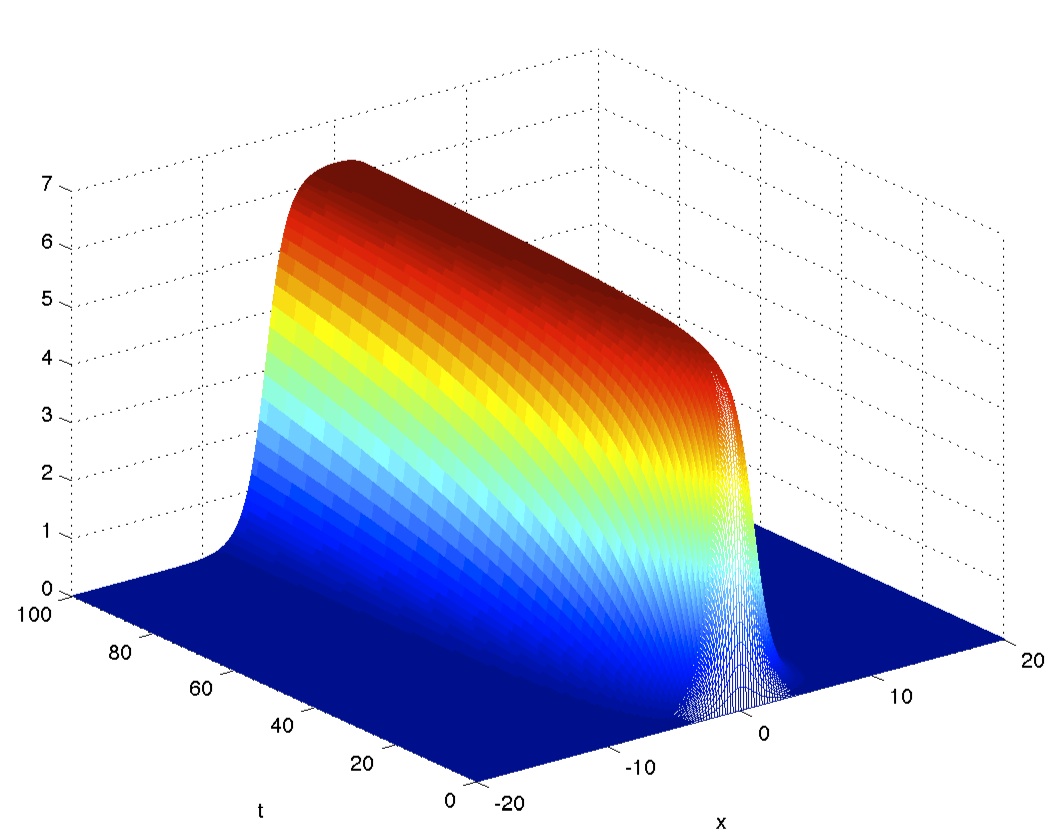}}
\caption{Soliton-like solution (\ref{sineGu})-(\ref{fi}) of problem (\ref{sineG})-(\ref{sineG0}), $\gamma=1$.}
\label{soli1}

\medskip
\centerline{\includegraphics[width=9cm,height=6cm]{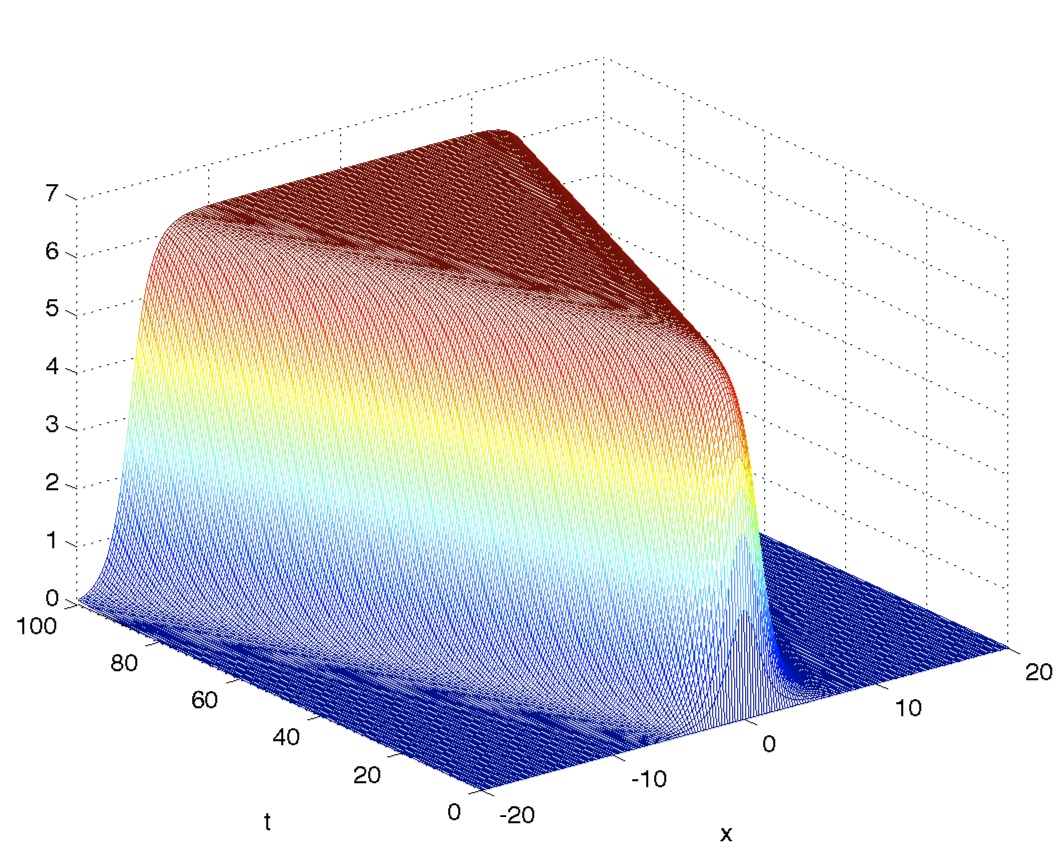}}
\caption{{\em Kink-antikink}, i.e., soliton-like solution (\ref{sineGu})-(\ref{fi}) of problem (\ref{sineG})-(\ref{sineG0}), $\gamma=0.99$.}
\label{soli2}
\end{figure}

We here consider a few numerical tests, concerning the so called {\em sine-Gordon} equation, which is in the form (\ref{wave}):
\begin{equation}\label{sineG}
u_{tt}(x,t) =u_{xx}(x,t)-\sin(u(x,t)).
\end{equation}
In particular, we shall consider {\em soliton-like} solutions, as described in \cite{W07}, defined by the initial conditions:
\begin{equation}\label{sineG0}
u(x,0) \equiv 0, \qquad u_t(x,0) = \frac{4}\gamma \sech\left(\frac{x}\gamma\right), \qquad \gamma>0.
\end{equation}
Depending on the value of the positive parameter $\gamma$, the solution is known to be given by:
\begin{equation}\label{sineGu}
u(x,t) = 4\atan\left[ \varphi(t;\gamma) \sech\left(\frac{x}\gamma\right)\right], 
\end{equation}
with
\begin{equation}\label{fi}
\varphi(t;\gamma) = \left\{ \begin{array}{ccc}
\frac{1}{\sqrt{\gamma^2-1}} \sin\left( \frac{\sqrt{\gamma^2-1}}\gamma t\right), &\qquad &\gamma>1,\\[5mm]
 t, & &\gamma=1,\\[5mm]
\frac{1}{\sqrt{1-\gamma^2}} \sinh\left( \frac{\sqrt{1-\gamma^2}}\gamma t\right), &\qquad &0<\gamma<1.\\
\end{array}\right.
\end{equation}
The three cases are shown in Figures~\ref{soli0}--\ref{soli2}, respectively: the first soliton (obtained for $\gamma>1$) is named {\em breather}, whereas the third one (obtained for $0<\gamma<1$) is named {\em kink-antikink}. Clearly, the case $\gamma=1$, shown in Figure~\ref{soli1}, separates the two different types of dynamics.

\begin{figure}
\centerline{\includegraphics[width=12cm,height=9cm]{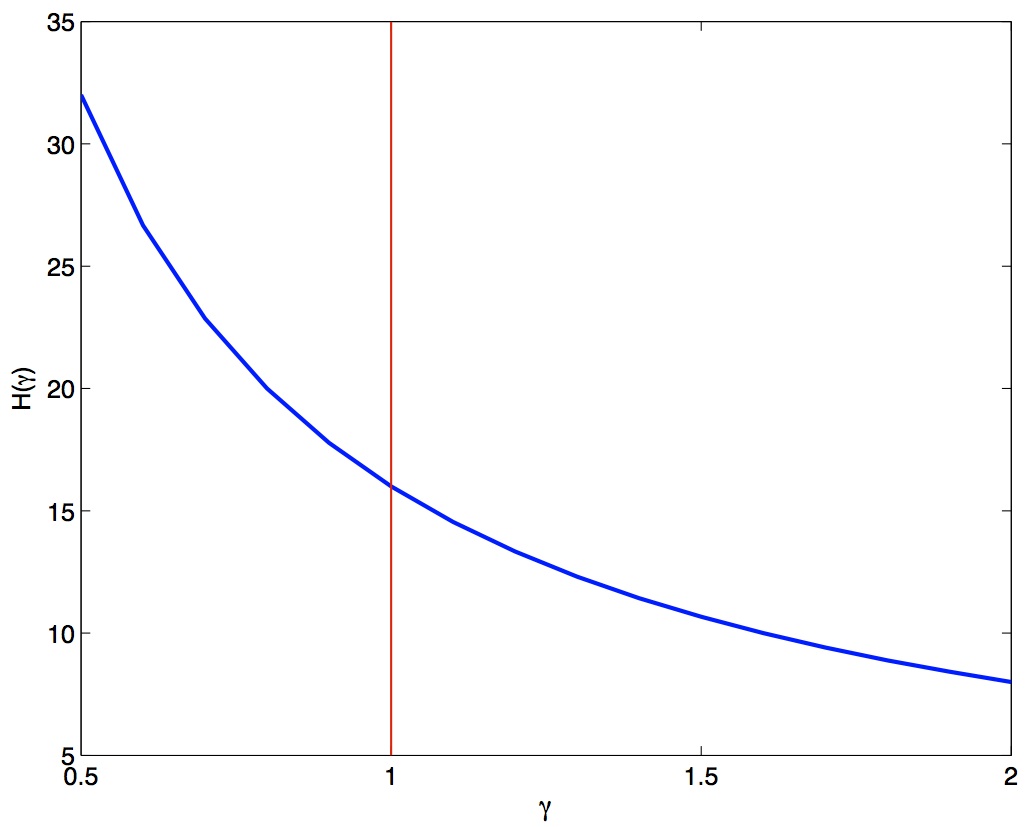}}
\caption{Hamiltonian for problem (\ref{sineG})-(\ref{sineG0}), as function of $\gamma$, for $x\in[-20,20]$.}
\label{hamil}

\bigskip
\centerline{\includegraphics[width=12cm,height=9cm]{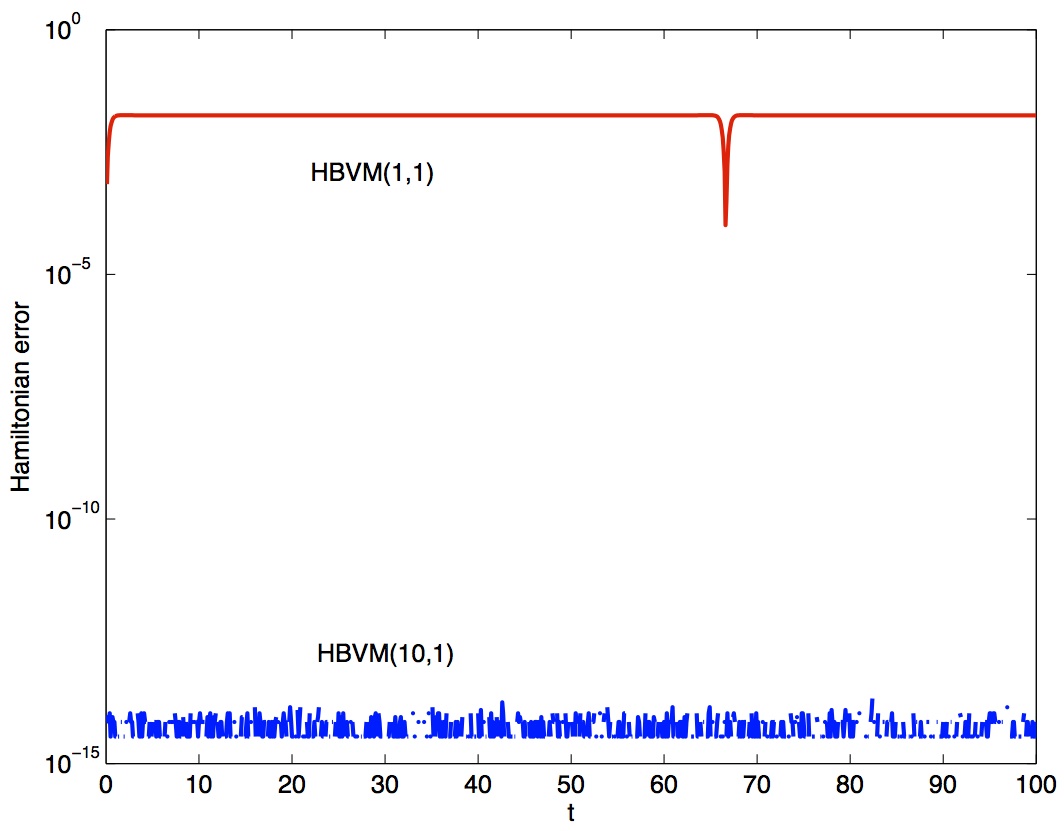}}
\caption{Hamiltonian error for the HBVM(1,1) and HBVM(10,1) methods, when solving problem (\ref{sineG})-(\ref{sineG0}), with $\gamma=1$, by using a stepsize $h=10^{-1}$.}
\label{Herr}
\end{figure}

Moreover, having fixed the space interval (we shall consider the interval $[-20,20]$), the Hamiltonian is a decreasing function of $\gamma$, as is shown in Figure~\ref{hamil}. This means that the value of the Hamiltonian (which is a constant of motion) characterizes the dynamics. Consequently, when $\gamma=1$, so that the Hamiltonian has a value $\simeq 16$, nearby values of the Hamiltonian will provide different types of soliton solutions. Consequently, energy conserving methods are expected to be useful, when numerically solving problem (\ref{sineG})-(\ref{sineG0}) with $\gamma=1$.

\begin{figure}[t]

\centerline{\includegraphics[width=12cm,height=9cm]{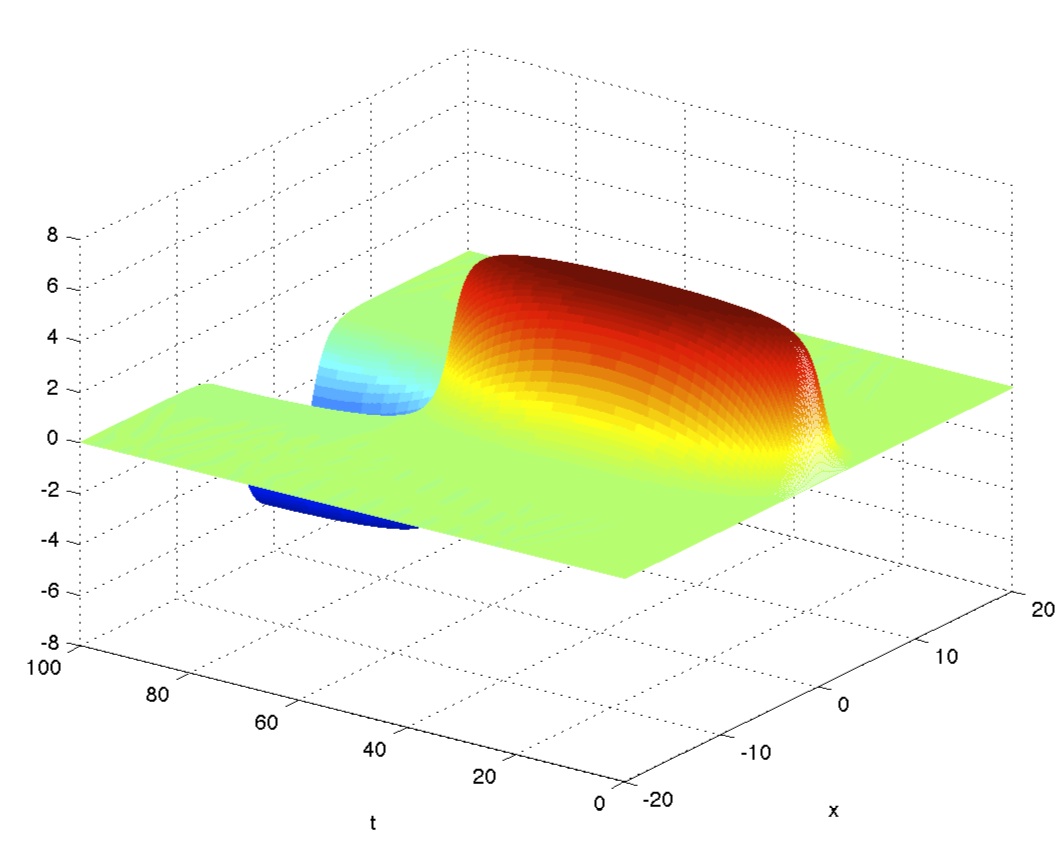}}
\caption{Numerical solution computed by the HBVM(1,1) method, when solving problem (\ref{sineG})-(\ref{sineG0}), with $\gamma=1$, by using a stepsize $h=10^{-1}$.}
\label{hbvm11}

\centerline{\includegraphics[width=12cm,height=9cm]{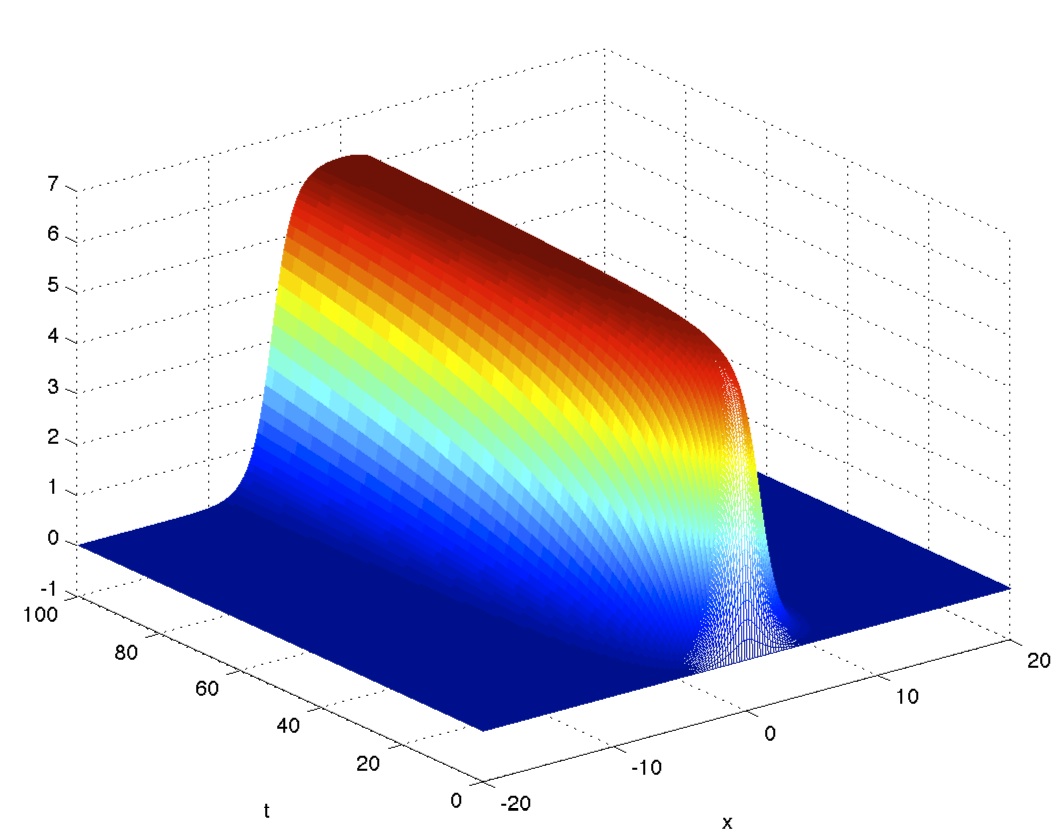}}
\caption{Numerical solution computed by the HBVM(10,1) method, when solving problem (\ref{sineG})-(\ref{sineG0}), with $\gamma=1$, by using a stepsize $h=10^{-1}$.}
\label{hbvm101}

\end{figure}

Let us then solve problem (\ref{sineG})-(\ref{sineG0}) with either periodic boundary conditions, or Dirichlet boundary conditions,  on the interval $[-20,20]$,\footnote{In fact, the resulting two discrete problems coincide, when using the approach in Section~\ref{firstapp} for coping with the case of Dirichlet boundary conditions.} by using:
\begin{itemize}
\item a trigonometric polynomial approximation of degree $N=100$;

\item $m=200$ equispaced mesh points in the given interval.
\end{itemize}
In so doing, the error in the initial condition is $\simeq 1.6\cdot 10^{-11}$, so that the initial profile is quite well matched. For the time integration, let us consider the following second-order methods, used with stepsize $h=10^{-1}$ for $10^3$ integration steps:
\begin{itemize}
\item the (symplectic) implicit mid-point rule, i.e., HBVM(1,1), for which the Hamiltonian error is $\simeq 1.8\cdot 10^{-2}$ (though without a drift);

\item the (practically) energy-conserving HBVM(10,1) method, for which the Hamiltonian error is $\simeq 2.1\cdot 10^{-14}$.
\end{itemize} 
The error in the numerical Hamiltonian is plotted in Figure~\ref{Herr}. In Figures~\ref{hbvm11} and \ref{hbvm101} we plot the numerical approximations to the solution computed by the HBVM(1,1) and HBVM(10,1) methods, respectively. As is clear, the former approximation is wrong, since the method has provided a {\em breather}-like solution, whereas the latter one well matches the continuous one (the maximum absolute error is $\simeq 4.6\cdot 10^{-3}$), thus confirming that energy conservation is an important issue, for such a problem.

\section{Conclusions}\label{final}

In this paper, we have studied the numerical solution of the nonlinear wave equation by using a Fourier discretization in space, also deriving a corresponding Hamiltonian formulation of the equation. Truncation of the Fourier expansion then leads to a corresponding truncated Hamiltonian, which turns out to be autonomous (thus conserved), when the problem is coupled with periodic bolundary-conditions. In case of different boundary conditions, the original approach can be suitably modified in order to recover a corresponding Hamiltonian problem with autonomous Hamiltonian. Energy-conserving methods in the HBVMs class can then be conveniently used for numerically solving the truncated problems. Energy-conservation turns out to be an interesting feature, for particular problems, possessing a soliton-like solution, as is confirmed by a few numerical tests.


\begin{thebibliography}{99}

\bibitem{AmSg05} P.\,Amodio, I.\,Sgura. High-order finite difference schemes for the solution of second-order BVPs. {\em Jour. Comput. Appl. Math.} {\bf 176} (2005) 59--76.

\bibitem{B01} J.P.\,Boyd. {\em  Chebyshev and Fourier spectral methods. Second edition.} Dover Publications, Inc., Mineola, NY, 2001.

\bibitem{BFI13} L.\,Brugnano, G.\,Frasca Caccia, F.\,Iavernaro. Efficient implementation of geometric integrators for separable Hamiltonian problems. {\em AIP Conf. Proc.} {\bf 1588} (2013) 734--737.

\bibitem{BFI14_0} L.\,Brugnano, G.\,Frasca Caccia, F.\,Iavernaro. Efficient implementation of Gauss collocation and Hamiltonian Boundary Value Methods. {\em Numer. Algor.} {\bf 65} (2014) 633--650.

\bibitem{BFI14} L.\,Brugnano, G.\,Frasca Caccia, F.\,Iavernaro. Energy conservation issues in the numerical solution of the nonlinear wave equation. {\em (submitted)}.

\bibitem{BIT09}
L.\,Brugnano, F.\,Iavernaro, D.\,Trigiante. Analysis of Hamiltonian Boundary Value Methods (HBVMs): a class of energy-preserving Runge-Kutta methods for the numerical solution of polynomial Hamiltonian systems.{\em Communications in Nonlinear Science and Numerical Simulation}
(2014), doi:\,\url{http://dx.doi.org/10.1016/j.cnsns.2014.05.030} \quad  (see also: \url{arXiv:0909.5659v1})

\bibitem{BIT09_1} L.\,Brugnano, F.\,Iavernaro, D.\,Trigiante. Hamiltonian BVMs (HBVMs): a family of  "drift-free" methods for integrating polynomial Hamiltonian systems. {\em AIP Conf. Proc.} {\bf 1168} (2009) 715--718.

\bibitem{BIT10}
L.\,Brugnano, F.\,Iavernaro, D.\,Trigiante. Hamiltonian Boundary Value Methods (Energy Preserving Discrete Line Methods). {\em Journal of Numerical Analysis, Industrial and Applied Mathematics} {\bf 5},1-2 (2010) 17--37.

\bibitem{BIT11}
L.\,Brugnano, F.\,Iavernaro, D.\,Trigiante. A note on the efficient implementation of Hamiltonian BVMs. {\em Journal of Computational and Applied Mathematics} {\bf 236} (2011) 375--383.

\bibitem{BIT12}
L.\,Brugnano, F.\,Iavernaro, D.\,Trigiante.  The Lack of Continuity and the Role of Infinite and Infinitesimal in Numerical Methods for ODEs: the Case of Symplecticity. {\em Applied Mathematics and Computation} {\bf 218} (2012) 8053--8063.

\bibitem{BIT12_1}
L.\,Brugnano, F.\,Iavernaro, D.\,Trigiante.  A simple framework for the derivation and analysis of effective one-step methods for ODEs. {\em Applied Mathematics and Computation} {\bf 218} (2012) 8475--8485.

\bibitem{BTbook} L.\,Brugnano, D.\,Trigiante. {\em Solving Differential Problems by Multistep Initial and Boundary Value Methods}, Gordon and Breach Science Publ., Amsterdam, 1998.

\bibitem{BrRe01} 
T.J.\,Bridges, S.\,Reich. Multi-symplectic integrators: numerical schemes for Hamiltonian PDEs that conserve symplecticity. {\em Physics Letters A} {\bf 284} (2001) 184--193.

\bibitem{CHQZ88} C.\,Canuto, M.Y.\,Hussaini, A.\,Quarteroni, T.A.\,Zang. {\em Spectral Methods in Fluid Dynamics}. Springer-Verlag, New York, 1988.

\bibitem{CHL08} D.\,Cohen, E.\,Hairer, C.\,Lubich. Conservation of energy, momentum and actions in numerical discretizations of non-linear wave equations. {\em Numer. Math.} {\bf 110}, no.\,2 (2008) 113--143.

\bibitem{DaBi08} G.\,Dahlquist, \AA.\,Bij\"{o}rk. {\em Numerical Methods in Scientific Computing, Vol.\,1}. SIAM, Philadelphia, 2008.

\bibitem{EvWe99} G.A.\,Evans, J.R.\,Webster. A comparison of some methods for the evaluation of highly oscillatory integrals. {\em Jour. Comput. Appl. Math.} {\bf 112} (1999) 55--69.

\bibitem{F12}
E.\,Faou. {\em Geometric Numerical Integration and Schr\"{o}dinger Equations.} Zurich, Switzerland: European Mathematical Society, 2012.

\bibitem{FMR06} 
J.\,Frank, B.E.\,Moore, S.\,Reich. Linear PDEs and Numerical Methods that Preserve a Multisymplectic Conservation Law. {\em SIAM J. Sci. Comput.} {\bf 28} (2006) 260--277.

\bibitem{FW78} 
B.\,Fornberg, G.B.\,Whitham.  A Numerical and Theoretical Study of Certain Nonlinear Wave Phenomena. {\em Proc. R. Soc. Lond. A} {\bf 289} (1978) 373--403.

\bibitem{HaWa96} 
E.\,Hairer, G.\,Wanner. {\em Solving Ordinary Differential Equations II. Stiff and Differential-Algebraic
Problems, 2nd edn}. Springer-Verlag, Berlin (1996)

\bibitem{KuRa09} A.\,Kurganov, J.\,Rauch. The Order of Accuracy of Quadrature Formulae for Periodic Functions.  {\em Advances in Phase Space Analysis of Partial Differential Equations}, A.\,Bove et al. (eds.),  Birkh\"auser, Boston,  2009.

\bibitem{IP07} F.\,Iavernaro, B.\,Pace. $s$-Stage Trapezoidal Methods for the Conservation of Hamiltonian Functions of Polynomial Type. {\em AIP Conf. Proc.} {\bf 936} (2007) 603--606.

\bibitem{IP08} F.\,Iavernaro, B.\,Pace. Conservative Block-Boundary Value Methods for the Solution of Polynomial Hamiltonian Systems. {\em AIP Conf. Proc.} {\bf 1048} (2008) 888--891.

\bibitem{IT09} F.\,Iavernaro, D.\,Trigiante. High-order symmetric schemes for the energy conservation of polynomial Hamiltonian problems. {\em Journal of Numerical Analysis, Industrial and Applied Mathematics} {\bf 4},1-2 (2009) 87--101.

\bibitem{IsSc04}
A.L.\,Islas, C.M.\,Schober. On the preservation of phase space structure under multisymplectic discretization. {\em Journal of Computational Physics} {\bf 197} (no. 2) (2004) 585--609.

\bibitem{LXW11} Z.\,Lv, M.\,Xue, Y.\,Wang. Legendre polynomials spectral approximation for the infinite-dimensional Hamiltonian systems. Math. Probl. in Engineering (2011) Article ID 824167, 13 pages.

\bibitem{Shen94} J.\,Shen. Efficient spectral-Galerkin method I. Direct solvers of second and fourth-order equations using Legendre polynomials. {\em SIAM Journal on Scientific Computing, 15(6),} 1489-1505

\bibitem{Shen95} J.\,Shen. Efficient spectral-Galerkin method II. Direct solvers of second and fourth-order equations using Chebyshev polynomials. {\em SIAM Journal on Scientific Computing, 16(1),} 74-87.

\bibitem{Wang91} D.\,Wang. Semi-discrete Fourier spectral approximations of infinite dimensional Hamiltonian systems and conservation laws. {\em Computers Math. Appl.} {\bf 21}, No.\,4 (1991) 63--75. 

\bibitem{WMGSS91} 
S.B.\,Wineberg, J.F.\,Mc\/Grath, E.F.\,Gabl, L.R.\,Scott,   C.E.\,Southwell. Implicit spectral methods for wave propogation problems. {\em J. Comp. Physics} {\bf 97} (1991) 311--336.

\bibitem{W07}
T.H.\,Wlodarczyk. {\em Stability and preservation properties of multisymplectic integrators}. PhD thesis, Department of Mathematics in the College of Sciences at the University of Central Florida, Orlando, Florida, 2007.
(\url{http://etd.fcla.edu/CF/CFE0001817/Wlodarczyk_Tomasz_H_200708_PhD.pdf})

\end{thebibliography}
 \end{document}